\documentclass[12pt]{amsart}
\usepackage[utf8]{inputenc}
\usepackage[T1]{fontenc}
\usepackage{amssymb,amsmath,amsthm}
\usepackage[a4paper,margin=1.9cm,top=2.5cm,bottom=2.5cm,centering,vcentering]{geometry}
\usepackage{color}
\definecolor{ao(english)}{rgb}{0.0, 0.0, 1.0}
\usepackage{hyperref}
\hypersetup{colorlinks=true, linkcolor=ao(english),citecolor=ao(english)}
\usepackage{orcidlink}
\usepackage{xcolor}

\definecolor{equationblue}{HTML}{244A73}
\definecolor{citationgreen}{HTML}{276347}
\definecolor{urlteal}{HTML}{176B75}

\hypersetup{
    colorlinks=true,
    linkcolor=blue,  
    citecolor=green!80!black, 
    urlcolor=blue    
}

\usepackage{breqn}

\allowdisplaybreaks
\numberwithin{equation}{section}

\newtheorem{theorem}{Theorem}[section]

\newtheorem{lemma}[theorem]{Lemma}
\newtheorem{remark}[theorem]{Remark}

\title[On Certain Congruences for Partition $k$-Tuples with $5$-Cores]{On Certain Congruences for Partition $k$-Tuples with $5$-Cores}

\author[S. Maity]{Saikat Maity \orcidlink{0009-0004-7188-471X}}
\address[S. Maity]{Mathematical and Physical Sciences division, School of Arts and Sciences, Ahmedabad University, Ahmedabad 380009, Gujarat, India}
\email{saikat.maths1729@gmail.com \& saikat.maity@ahduni.edu.in}

\begin{document}

\begin{abstract}
Let $A_{t,k}(n)$ denote the number of partition $k$-tuples of $n$ with $t$-cores. Using fifth-order modular equations and the Atkin $U_5$-operator, we prove the congruences
$A_{5,125i+21}(25n+r)\equiv 0 \pmod{5^3}$ for $r\in\{9,14,19,24\}$,
$A_{5,125i+41}(25n+r)\equiv 0 \pmod{5^3}$ for $r\in\{19,24\}$, and
$A_{5,125i+111}(125n+r)\equiv 0 \pmod{5^4}$ for $r\in\{39,64,89,114\}$,
which were stated by Tang \cite{tang2025iterative}.
\end{abstract}

\subjclass[2020]{11P83; 05A17}
\keywords{Partitions, $t$-cores, partition $k$-tuples, congruences, modular equations}

\maketitle

\section{Introduction}

A \emph{partition} of a positive integer $n$ is a finite non-increasing sequence of positive integers
$\lambda = (\lambda_1, \lambda_2, \dots, \lambda_\ell)$ with $\sum_{i=1}^{\ell} \lambda_i = n$.
Given a partition $\lambda$, its Ferrers--Young diagram is the left-justified array of dots with $\lambda_i$ dots in the $i$-th row. The \emph{hook number} of a dot is the total number of dots directly to the right and directly below the dot, including the dot itself. A partition is called a \emph{$t$-core partition} for $t\ge 2$ if none of its hook numbers is divisible by $t$. A \emph{partition $k$-tuple} of $n$ is an ordered $k$-tuple of partitions
$(\Lambda_1, \Lambda_2, \dots, \Lambda_k)$
such that
$|\Lambda_1| + |\Lambda_2| + \cdots + |\Lambda_k| = n$.
If each $\Lambda_i$ is a $t$-core partition, the tuple is called a \emph{$t$-core partition $k$-tuple}.

Let $A_{t,k}(n)$ denote the number of partition $k$-tuples of $n$ with $t$-cores. Garvan, Kim, and Stanton \cite{garvan1990cranks} showed that
\begin{equation} \label{genfunc}
\sum_{n=0}^{\infty} A_{t,k}(n)q^n = \frac{(q^t;q^t)_{\infty}^{kt}}{(q;q)_{\infty}^k},
\end{equation}
where $(a;q)_\infty := \prod_{i=0}^{\infty} (1-aq^i)$ and $|q|<1$.

Congruences for $A_{t,k}(n)$ have been studied in a number of cases. For $k=1$, see Hirschhorn and Sellers \cite{HirschhornSellers1996,HirschhornSellers1999} and Radu and Sellers \cite{RaduSellers2011}; for $t=3$, see Lin \cite{Lin2014}, Xia \cite{Xia2015} and Wang \cite{Wang2016a,Wang2016b}. For $t=5$, Ranganatha \cite{Ranganatha2016} and Majid and Fathima \cite{MajidFathima2022} proved congruences for $k=2$ and $k=3$, respectively. Saikia, Sarma and Talukdar \cite{SaikiaSarmaTalukdar2024} treated $k=2,3,4$, obtained families of congruences for $k$-tuples with $p$-cores ($p$ prime), and stated further conjectural congruences. See also \cite{Ranganatha2021} for general properties of $k$-tuple $t$-core partitions.

Tang \cite{tang2025iterative} introduced an iterative approach, based on fifth-order modular equations, for proving congruences modulo powers of $5$ for $A_{5,k}(n)$. He proved six of the conjectures in \cite{SaikiaSarmaTalukdar2024} and stated several further congruences, among them (4.4), (4.5) and (4.6), without proof. 
Liang and Tang \cite{Liang_Tang_2025} used a recurrence coming from a fifth-order modular equation to obtain many families of congruences for $4\le k\le 25$, and Ananya, Keerthana and Ranganatha \cite{AnanyaKeerthanaRanganatha2026} recently proved congruences for $A_{5,5^ri+j}(n)$ with $1\le j\le 8$.

The purpose of this paper is to prove (4.4), (4.5) and (4.6) of \cite{tang2025iterative} (Theorems~\ref{thm:cong4.4}--\ref{thm:cong4.6}). The proofs use the Atkin $U_5$-operator together with a recurrence for $U_5(\zeta^k)$ coming from a fifth-order modular equation, as in \cite{Liang_Tang_2025}.

Section~\ref{prelim} collects the required $q$-series facts and the values of $U_5(\zeta^k)$ that are needed. Section~\ref{mainresults} contains the proofs.

\section{Preliminaries}\label{prelim}
Throughout, $|q|<1$ and
\[
E(q) := (q;q)_\infty = \prod_{n=1}^{\infty}(1-q^n).
\]
All congruences between power series are understood coefficientwise.

\begin{lemma}\label{lem:2.1}
\textbf{(Euler's Pentagonal Number Theorem)}
\begin{equation}
\label{eq:2.1}
E(q)=\sum_{n=-\infty}^{\infty}(-1)^n q^{n(3n-1)/2}.
\end{equation}
\end{lemma}

\begin{proof}
See \cite[Corollary 1.7]{Andrews1976}.
\end{proof}

For a Laurent series $\sum_{n} a_n q^n$, the Atkin $U_k$-operator is defined by
\[
U_k\left(\sum_{n=-\infty}^{\infty} a_n q^n\right) = \sum_{n=-\infty}^{\infty} a_{kn} q^n.
\]
Note that $U_k\big(f(q^k)\,g(q)\big)=f(q)\,U_k\big(g(q)\big)$.

\begin{lemma}\label{lem:2.2}
Let
\[
\zeta = \frac{E(q)}{qE(q^{25})},
\qquad
\xi = \frac{E(q)^6}{qE(q^5)^6}.
\]
Then for any integer $i\ge0$, we have
\begin{align}
U_5(\zeta^{i+5})
&=
-5U_5(\zeta^{i+4})-15U_5(\zeta^{i+3})-25U_5(\zeta^{i+2}) -25U_5(\zeta^{i+1})+\xi U_5(\zeta^{i}),
\label{eq:2.3}
\end{align}
where
\begin{equation}
\label{eq:2.4}
U_5(\zeta^0)=1,\quad
U_5(\zeta)=-1,\quad
U_5(\zeta^2)=-1,\quad
U_5(\zeta^3)=5,\quad
U_5(\zeta^4)=-5.
\end{equation}
\end{lemma}

\begin{proof}
The initial condition $U_5(\zeta^0) = 1$ is trivial. The remaining four initial values in \eqref{eq:2.4} can be computed using \cite[Eq. (6.3.3)]{Hirschhorn2017}. 
Moreover, the following fifth-order modular equation holds
\cite[Eq.~(6.2.14)]{Hirschhorn2017}:
\begin{equation}\label{eq:2.5}
\zeta^5 = -5\zeta^4-15\zeta^3-25\zeta^2-25\zeta+\frac{E(q^5)^6}{q^5E(q^{25})^6}.
\end{equation}
Multiplying \eqref{eq:2.5} by $\zeta^i$ and applying the $U_5$-operator
yields \eqref{eq:2.3}, since the last term of \eqref{eq:2.5} equals $\xi(q^5)$ and $U_5(\xi(q^5)\zeta^i)=\xi(q)\,U_5(\zeta^i)$.
\end{proof}

\begin{remark}\label{rem:2.5}
It follows from Lemma~\ref{lem:2.2} by induction that, for each integer $k \ge 1$,
$U_5(\zeta^k)=P_k(\xi)$,
where $P_k$ is a polynomial with integer coefficients of degree at most $\lfloor k/5 \rfloor$.
\end{remark}

\begin{lemma}\label{lem:binomial}
Let $p$ be a prime and $k$ a positive integer. Then
\[
E(q)^{p^k}
\equiv
E(q^p)^{p^{k-1}}
\pmod{p^k}.
\]
\end{lemma}

\begin{proof}
We have $(1-x)^p \equiv 1-x^p \pmod{p}$. Since $A\equiv B \pmod{p^j}$ implies $A^p\equiv B^p \pmod{p^{j+1}}$, induction on $k$ gives
\[
(1-x)^{p^k}\equiv (1-x^p)^{p^{k-1}} \pmod{p^k}.
\]
Taking $x=q^n$ and the product over $n\ge1$ completes the proof.
\end{proof}

Since $E(q)$ and $E(q^p)$ have constant term $1$, Lemma~\ref{lem:binomial} remains valid after raising both sides to any integer power, positive or negative.

Using the recurrence \eqref{eq:2.3} and the initial conditions \eqref{eq:2.4}, we compute
\begin{align}
U_5(\zeta^{84}) &= 1647135\xi^{16} + 625\sum_{k=0}^{15} c_{84,k}\xi^k, \label{eq:u5zeta84}\\
U_5(\zeta^{104}) &= 4057170\xi^{20} + 3125\sum_{k=0}^{19} c_{104,k}\xi^k, \label{eq:u5zeta104}\\
U_5(\zeta^{514}) &= 2807114625\xi^{102} + 15625\sum_{k=0}^{101} c_{514,k}\xi^k, \label{eq:u5zeta514}
\end{align}
where the $c_{i,k}$ are integers. For the proof of Theorem~\ref{thm:cong4.6} with arbitrary $i$ we also need the values at $14$, $139$, $264$, $389$:
\begin{align}
U_5(\zeta^{14}) &= -15625, \label{eq:u5zeta14}\\
U_5(\zeta^{139}) &= 13587250\,\xi^{27} + 15625\sum_{k=0}^{26} c_{139,k}\xi^k, \label{eq:u5zeta139}\\
U_5(\zeta^{264}) &= 188768250\,\xi^{52} + 15625\sum_{k=0}^{51} c_{264,k}\xi^k, \label{eq:u5zeta264}\\
U_5(\zeta^{389}) &= 910308625\,\xi^{77} + 15625\sum_{k=0}^{76} c_{389,k}\xi^k. \label{eq:u5zeta389}
\end{align}
Reducing \eqref{eq:u5zeta14}--\eqref{eq:u5zeta389} and \eqref{eq:u5zeta514} modulo $5^4$ (note $13587250\equiv 375$, $188768250\equiv 125$, $910308625\equiv 500$, $2807114625\equiv 250 \pmod{5^4}$ and $5^4\mid 15625$), we obtain, for $j=0,1,2,3,4$,
\begin{equation}\label{eq:u5zeta14j}
U_5(\zeta^{14+125j}) \equiv 375\,j\,\xi^{25j+2} \pmod{5^4}.
\end{equation}

\section{Main Results} \label{mainresults}

\begin{theorem}[Congruence (4.4), \cite{tang2025iterative}]
\label{thm:cong4.4}
For all $n \geq 0$ and $i \geq 0$,
\begin{equation}
\label{eq:cong4.4}
A_{5,125i+21}(25n+r) \equiv 0 \pmod{5^3}, \qquad r\in\{9,14,19,24\}.
\end{equation}
\end{theorem}

\begin{proof}
By \eqref{genfunc},
\begin{equation*}
\sum_{n\geq 0} A_{5,k}(n)q^n = \frac{E(q^5)^{5k}}{E(q)^k}.
\end{equation*}
Setting $k = 125i + 21$, we obtain
\begin{align*}
\sum_{n \geq 0} A_{5,125i+21}(n) q^{n}
&= \frac{E(q^{5})^{625i+105}}{E(q)^{125i}E(q)^{21}} \\
&\equiv \frac{E(q^{5})^{625i+105}}{E(q^{5})^{25i}E(q)^{21}}
\pmod{5^{3}}
\quad \text{(by Lemma~\ref{lem:binomial})} \\
&= \frac{E(q^{5})^{600i+105}}{E(q)^{125}} E(q)^{104} \\
&\equiv \frac{E(q^{5})^{600i+105}}{E(q^{5})^{25}} E(q)^{104}
\pmod{5^{3}}
\quad \text{(by Lemma~\ref{lem:binomial})} \\
&= E(q^{5})^{600i+80} E(q)^{104}.
\end{align*}
Multiplying by $q^{-4}$ and using $\zeta = \frac{E(q)}{qE(q^{25})}$, we have
\begin{align*}
\sum_{n=0}^{\infty} A_{5,125i+21}(n)q^{n-4}
\equiv
q^{100}E(q^5)^{600i+80}E(q^{25})^{104}\zeta^{104}
\pmod{5^3}.
\end{align*}
Applying the $U_5$-operator to both sides gives
\begin{equation}\label{eq1:cong4.4}
\sum_{n=0}^{\infty} A_{5,125i+21}(5n+4)q^n \equiv q^{20}E(q)^{600i+80}E(q^5)^{104}U_5(\zeta^{104})\pmod{5^3}.
\end{equation}

By \eqref{eq:u5zeta104}, since $4057170\equiv 45 \pmod{5^3}$ and $5^3\mid 3125$,
\begin{equation*}
U_5(\zeta^{104}) \equiv 45\xi^{20} \pmod{5^3}.
\end{equation*}
Substituting this into \eqref{eq1:cong4.4} and then $\xi = \frac{E(q)^6}{qE(q^5)^6}$, we obtain
\begin{align*}
\sum_{n=0}^{\infty} A_{5,125i+21}(5n+4)q^n &\equiv 45q^{20}E(q)^{600i+80}E(q^5)^{104} \left(\frac{E(q)^6}{qE(q^5)^6}\right)^{20} \pmod{5^3}\\
&= 45E(q)^{600i+200}E(q^5)^{-16}.
\end{align*}
Since $5\mid 45$, it suffices to know $E(q)^{600i+200}$ modulo $5^2$. By Lemma~\ref{lem:binomial} with $k=2$, $E(q)^{25}\equiv E(q^5)^5 \pmod{5^2}$, so $E(q)^{600i+200}=\big(E(q)^{25}\big)^{24i+8}\equiv E(q^5)^{120i+40} \pmod{5^2}$. Hence
\begin{align*}
\sum_{n=0}^{\infty} A_{5,125i+21}(5n+4)q^n &\equiv 45E(q^5)^{120i+24} \pmod{5^3}.
\end{align*}
The right-hand side is a power series in $q^5$, so the coefficients of $q^{5n+s}$ vanish modulo $5^3$ for $s\in\{1,2,3,4\}$. Since $5(5n+s)+4=25n+5s+4$, this gives \eqref{eq:cong4.4}.
\end{proof}

\begin{theorem}[Congruence (4.5), \cite{tang2025iterative}]
\label{thm:cong4.5}
For all $n \geq 0$ and $i \geq 0$,
\begin{equation}
\label{eq:cong4.5}
A_{5,125i+41}(25n+r) \equiv 0 \pmod{5^3}, \qquad r\in\{19,24\}.
\end{equation}
\end{theorem}

\begin{proof}
Setting $k = 125i + 41$ in \eqref{genfunc}, we have
\begin{align*}
\sum_{n \geq 0} A_{5,125i+41}(n) q^{n}
&= \frac{E(q^{5})^{625i+205}}{E(q)^{125i}E(q)^{41}} \\
&\equiv \frac{E(q^{5})^{625i+205}}{E(q^{5})^{25i}E(q)^{41}}
\pmod{5^{3}}
\quad \text{(by Lemma~\ref{lem:binomial})} \\
&= \frac{E(q^{5})^{600i+205}}{E(q)^{125}} E(q)^{84} \\
&\equiv \frac{E(q^{5})^{600i+205}}{E(q^{5})^{25}} E(q)^{84}
\pmod{5^{3}}
\quad \text{(by Lemma~\ref{lem:binomial})} \\
&= E(q^{5})^{600i+180} E(q)^{84}.
\end{align*}
Multiplying by $q^{-4}$ and using $\zeta = \frac{E(q)}{qE(q^{25})}$, we have
\begin{align*}
\sum_{n=0}^{\infty} A_{5,125i+41}(n)q^{n-4} &\equiv q^{80} E(q^5)^{600i+180} E(q^{25})^{84} \zeta^{84} \pmod{5^3}.
\end{align*}
Applying the $U_5$-operator to both sides gives
\begin{equation}\label{eq1:cong4.5}
\sum_{n=0}^{\infty} A_{5,125i+41}(5n+4)q^n \equiv q^{16} E(q)^{600i+180}E(q^5)^{84}U_5(\zeta^{84}) \pmod{5^3}.
\end{equation}
By \eqref{eq:u5zeta84}, since $1647135\equiv 10 \pmod{5^3}$ and $5^3\mid 625$,
\begin{equation*}
U_5(\zeta^{84}) \equiv 10\xi^{16} \pmod{5^3}.
\end{equation*}
Substituting this into \eqref{eq1:cong4.5} and then $\xi = \frac{E(q)^6}{qE(q^5)^6}$, we obtain
\begin{align*}
\sum_{n=0}^{\infty} A_{5,125i+41}(5n+4)q^n &\equiv 10q^{16} E(q)^{600i+180}E(q^5)^{84} \left(\frac{E(q)^6}{q E(q^5)^6}\right)^{16}\\
&= 10 E(q)^{600i+276}E(q^5)^{-12}.
\end{align*}
Since $5\mid 10$, it suffices to know $E(q)^{600i+275}=\big(E(q)^{25}\big)^{24i+11}$ modulo $5^2$, and by Lemma~\ref{lem:binomial} with $k=2$ it is congruent to $E(q^5)^{120i+55}$. Hence
\begin{align*}
\sum_{n=0}^{\infty} A_{5,125i+41}(5n+4)q^n &\equiv 10 E(q^5)^{120i+43}E(q) \pmod{5^3}.
\end{align*}
By \eqref{eq:2.1}, the coefficient of $q^m$ in $E(q)$ vanishes unless $m=n(3n-1)/2$ for some integer $n$. Then $24m+1=(6n-1)^2$ is a square, and since $24m+1\equiv 1-m \pmod 5$ and the squares modulo $5$ are $0,1,4$, we get $m\equiv 0,1,2 \pmod 5$. Writing $\sum_{n\ge0}p_1(n)q^n=E(q)$, this means
\begin{equation} \label{p1}
p_1(5n+3)=p_1(5n+4)=0 \qquad (n\ge 0).
\end{equation}
Since $E(q^5)^{120i+43}$ is a series in $q^5$, \eqref{p1} shows that the coefficients of $q^{5n+3}$ and $q^{5n+4}$ in the right-hand side above vanish. As $5(5n+s)+4=25n+5s+4$, taking $s=3,4$ gives \eqref{eq:cong4.5}.
\end{proof}

\begin{theorem}[Congruence (4.6), \cite{tang2025iterative}]
\label{thm:cong4.6}
For all $n \geq 0$ and $i \geq 0$,
\begin{equation}
\label{eq:cong4.6}
A_{5,125i+111}(125n+r) \equiv 0 \pmod{5^4}, \qquad r\in\{39,64,89,114\}.
\end{equation}
\end{theorem}

\begin{proof}
By \eqref{genfunc}, for $k = 125i + 111$ we have
\[
\sum_{n \geq 0} A_{5,125i+111}(n) q^{n} = \frac{E(q^{5})^{625i+555}}{E(q)^{125i+111}}.
\]
Write $i=5a+b$ with $0\le b\le 4$, and put $j=4-b$ and $N=514-125b=14+125j$. Then $125i+111=625(a+1)-N$, so
\begin{align*}
\sum_{n \geq 0} A_{5,125i+111}(n) q^{n}
&= \frac{E(q^{5})^{625i+555}}{E(q)^{625(a+1)}}\,E(q)^{N} \\
&\equiv \frac{E(q^{5})^{625i+555}}{E(q^{5})^{125(a+1)}}\,E(q)^{N}
\pmod{5^{4}}
\quad \text{(by Lemma~\ref{lem:binomial} with $k=4$)} \\
&= E(q^{5})^{e}\,E(q)^{N},
\end{align*}
where $e=625i-125a+430=600i+25b+430$. (For $b=0$ this is the exponent $600i+430$ together with $E(q)^{514}$.) The exponent $625(a+1)$ is chosen so that Lemma~\ref{lem:binomial} is applied with $k=4$ directly, which gives the congruence modulo $5^4$ for every $i\ge0$.

Multiplying by $q^{-4}$ and using $E(q)^N=q^NE(q^{25})^N\zeta^N$, we obtain
\begin{align*}
\sum_{n=0}^{\infty} A_{5,125i+111}(n)q^{n-4} \equiv q^{N-4}E(q^5)^{e}E(q^{25})^{N}\zeta^{N} \pmod{5^4},
\end{align*}
and $N-4=5(25j+2)$. Applying the $U_5$-operator gives
\begin{equation}\label{eq1:cong4.6}
\sum_{n=0}^{\infty} A_{5,125i+111}(5n+4)q^n \equiv q^{25j+2}E(q)^{e}E(q^5)^{N}U_5(\zeta^{N}) \pmod{5^4}.
\end{equation}
By \eqref{eq:u5zeta14j}, $U_5(\zeta^{N})\equiv 375\,j\,\xi^{25j+2} \pmod{5^4}$. Substituting this and $\xi = \frac{E(q)^6}{qE(q^5)^6}$, we get
\begin{align*}
\sum_{n=0}^{\infty} A_{5,125i+111}(5n+4)q^n
&\equiv 375\,j\,q^{25j+2}E(q)^{e}E(q^5)^{N}
\left(\frac{E(q)^6}{qE(q^5)^6}\right)^{25j+2} \pmod{5^4}\\
&= 375\,j\,E(q)^{600i-125b+1042}E(q^5)^{2-25j},
\end{align*}
since $e+6(25j+2)=600i+25b+430+150j+12=600i-125b+1042$ (using $j=4-b$) and $N-6(25j+2)=2-25j$. For $b=0$ this is $250E(q)^{600i+1042}E(q^5)^{-98}$, as in the original computation.
Since $375=3\cdot 5^3$, it suffices to work modulo $5$ in the remaining factor. By Lemma~\ref{lem:binomial} with $k=1$, $E(q)^{600i-125b+1040}=\big(E(q)^{5}\big)^{120i-25b+208}\equiv E(q^5)^{120i-25b+208} \pmod 5$, so
\begin{align*}
\sum_{n=0}^{\infty} A_{5,125i+111}(5n+4)q^n
&\equiv 375\,j\,E(q^5)^{120i+110} E(q)^{2} \pmod{5^4}\\
&= 375\,j\,q^2 E(q^5)^{120i+110} E(q^{25})^{2} \zeta^{2},
\end{align*}
because $(120i-25b+208)+2-25j=120i+110$.
Hence
\begin{align*}
    \sum_{n=0}^{\infty} A_{5,125i+111}(5n+4)q^{n-2} \equiv 375\,j\,E(q^5)^{120i+110} E(q^{25})^{2} \zeta^{2} \pmod{5^4}.
\end{align*}
Applying the $U_5$-operator again, and using $U_5(\zeta^2)=-1$ from Lemma~\ref{lem:2.2}, we get
\begin{align*}
    \sum_{n=0}^{\infty} A_{5,125i+111}(25n+14)q^{n} &\equiv 375\,j\,E(q)^{120i+110} E(q^{5})^{2} U_5(\zeta^{2})\\
    &= - 375\,j\,E(q)^{120i+110} E(q^{5})^{2} \\
    &\equiv - 375\,j\,E(q^5)^{24i+22} E(q^{5})^{2} \pmod{5^{4}} \quad \text{(by Lemma~\ref{lem:binomial} with $k=1$)}\\
    &= - 375\,j\,E(q^5)^{24i+24}.
\end{align*}
The right-hand side is a power series in $q^5$, so the coefficients of $q^{5n+s}$ vanish modulo $5^4$ for $s\in\{1,2,3,4\}$. Since $25(5n+s)+14=125n+25s+14$, taking $s=1,2,3,4$ gives \eqref{eq:cong4.6}.
\end{proof}

\section*{Acknowledgements}
The author thanks Manjil P. Saikia for helpful discussions and comments.

\end{document}